\documentclass[12pt]{article}
\usepackage{amsmath}
\usepackage{amssymb}
\newcommand{\D}{{\cal D}}
\newcommand{\tauxf}{\|\tau_xf-f\|}
\newcommand{\hk}{{\cal HK}}
\newcommand{\bv}{{\cal BV}}
\newcommand{\R}{{\mathbb R}}

\newcommand{\N}{{\mathbb N}}
\newcommand{\fn}{\!:\!}
\newcommand{\lsup}{\sup\limits}
\newcommand{\Qed}{\mbox{~~~$\blacksquare$}}
\newtheorem{theorem}{Theorem}
\newtheorem{lemma}[theorem]{Lemma}

\newtheorem{corollary}[theorem]{Corollary}

\newtheorem{example}[theorem]{Example}
\begin{document}
\begin{center}
{\large\bf Continuity in the Alexiewicz norm}
\vskip.25in
Erik Talvila\footnote{Supported by the
Natural Sciences and Engineering Research Council of Canada.
}\\ 
{\footnotesize
Department of Mathematics and Statistics \\
University College of the Fraser Valley\\
Abbotsford, BC Canada V2S 7M8\\
Erik.Talvila@ucfv.ca\\

{\it Dedicated to Jaroslav Kurzweil on his $80^{th}$ birthday.}\\
}
\end{center}
{\footnotesize
\noindent
{\bf Abstract.} If $f$ is a Henstock--Kurzweil integrable function on the
real line, the Alexiewicz
norm of $f$ is $\|f\|=\sup_I|\int_I f|$ where the supremum is taken over
all intervals $I\subset\R$.
Define the translation $\tau_x$ by $\tau_xf(y)=f(y-x)$.  Then $\|\tau_xf-f\|$ tends to
$0$ as $x$ tends to $0$, i.e., $f$ is continuous in the Alexiewicz norm.
For particular
functions,  $\|\tau_xf-f\|$ can tend to 0 arbitrarily slowly.
In general,
$\|\tau_xf-f\|\geq {\rm osc}f\,|x|$ as $x\to 0$, where ${\rm osc}f$ is the oscillation
of $f$.  It is shown
that if $F$ is a primitive of $f$ then  $\|\tau_xF-F\|\leq \|f\||x|$.  
An example shows that the function $y\mapsto \tau_xF(y)-F(y)$ need not
be in $L^1$.  However, if $f\in L^1$ then $\|\tau_xF-F\|_1\leq \|f\|_1|x|$.
For a positive weight function $w$ on the real line, necessary and sufficient
conditions on $w$ are given so that $\|(\tau_xf-f)w\|\to 0$ as $x\to 0$ whenever
$fw$ is Henstock--Kurzweil integrable.  Applications are made to the
Poisson integral on the disc and half-plane.  
All of the results also
hold with the distributional Denjoy integral, which arises from
the completion of the
space of Henstock--Kurzweil integrable functions as a subspace of Schwartz
distributions.
}\\

\section{Introduction}
For $f\fn\R\to\R$ define the translation by $\tau_xf(y)=f(y-x)$ for
$x,y\in\R$.  If $f\in L^p$ ($1\leq p<\infty$) then it is a well
known result of Lebesgue integration that $f$ is continuous in the
$p$-norm, i.e., $\lim_{x\to 0}\|\tau_xf-f\|_p=0$.  For example, see
\cite[Lemma~6.3.5]{stroock}.
In this paper we consider continuity of Henstock--Kurzweil
integrable functions in Alexiewicz and weighted Alexiewicz norms on
the real line.  Let $\hk$ be the set of functions $f\fn\R\to\R$ that are Henstock--Kurzweil
integrable.  The Alexiewicz norm of $f\in\hk$ is defined
$\|f\|=\sup_I|\int_If|$ where the supremum is over all intervals 
$I\subset\R$.
Identifying functions almost everywhere, $\hk$ becomes
a normed linear space under $\|\cdot\|$ that is barrelled
but not complete.
See \cite{lee} and \cite{swartz}
for a discussion of the Henstock--Kurzweil
integral and the Alexiewicz norm.
It is shown below that translations are continuous in norm and
that for $f\in\hk$ we have
$\|\tau_xf-f\|\geq {\rm osc} f\,|x|$ where ${\rm osc} f$ is the oscillation
of $f$.  For particular $f\in\hk$ the
quantity  $\|\tau_xf-f\|$ can tend to 0 arbitrarily slowly.  If $F$ is a primitive of
$f$ then $\|\tau_xF-F\|\leq \|f\||x|$. 
An example shows
that if $f\in \hk$ then the function defined by $y\mapsto \tau_xF(y)-F(y)$
need not be in $L^1$ but if $f\in L^1$ then $\|\tau_xF-F\|_1\leq \|f\|_1|x|$. 
For a positive weight function $w$ on the real line, necessary and sufficient
conditions on $w$ are given so that $\|(\tau_xf-f)w\|\to 0$ as $x\to 0$ whenever
$fw$ is Henstock--Kurzweil integrable. The necessary and sufficient conditions
involve properties of the function $g_x(y)=w(y+x)/w(y)$.  Sufficient conditions
are given on $w$ for $\|(\tau_xf-f)w\|\to 0$.
Applications to the Dirichlet problem in the disc and half-plane are
given.  

All of the results also hold when we use the distributional
Denjoy integral.  Define ${\cal A}$ to be the completion of $\hk$ with
respect to $\|\cdot\|$.
Then ${\cal A}$ is a subspace of the space of
Schwartz distributions.  Distribution $f$ is in ${\cal A}$ if there is
function $F$ continuous on the extended real line such that $F'=f$ as
a distributional derivative.  
For details on this integral see \cite{talviladenjoy}.

First we prove continuity in the Alexiewicz norm.

\begin{theorem}\label{theorem1}
Let $f\in\hk$.  For $x,y\in\R$ define
$\tau_xf(y)=f(y-x)$.  Then $\|\tau_xf-f\|\to 0$ as $x\to 0$.
\end{theorem}
\noindent
{\bf Proof:} Let $x,\alpha,\beta\in\R$.  Then
$
\int_{\alpha}^{\beta}(\tau_xf-f) = \int_{\alpha-x}^{\beta-x}\!
f-\int_\alpha^{\beta}\!f$. Write $F(x)=\int_{-\infty}^xf$.
Taking the supremum over $\alpha$ and $\beta$,
\begin{eqnarray*}
\|\tau_xf-f\| & \leq & \lsup_{\beta\in\R}|F(\beta -x)-F(\beta)|+
\lsup_{\alpha\in\R}|F(\alpha -x)-F(\alpha)|\\
 & \to & 0 \text{ as } x\to 0 \text{ since } F \text{ is uniformly
continuous on } \R.
\Qed
\end{eqnarray*}

Notice that for each $x\in\R$, the translation $\tau_x$ is an
isometry on $\hk$, i.e., it is a homeomorphism such that
$\|\tau_xf\|=\|f\|$.  It is also clear
that we have continuity at each point: for each $x_0\in \R$, $\|\tau_xf
-\tau_{x_0}f\|\to 0$ as $x\to x_0$.

The theorem also applies on any interval $I\subset\R$.
Restrict $\alpha$ and $\beta$ to lie in $I$ and extend $f$ to be $0$
outside $I$.  Or, one could use a periodic extension.  The same results
also hold for the equivalent norm $\|f\|=\sup_{x\in\R}|\int_{-\infty}^x f|$.

Under the Alexiewicz norm, the space of
Henstock--Kurzweil integrable functions  is not
complete.  Its completion with respect to the norm
$\|f\|=\sup_{x\in\R}|\int_{-\infty}^x f|$ is the subspace of distributions that
are the distributional derivative of a function in 
$\tilde{C}:=\{F\fn\R\to\R\mid F\in C^0(\R), \lim_{x\to{-\infty}}F(x)=0,
\lim_{x\to{\infty}}F(x)\in\R\}$, i.e., they are 
distributions of order 1. See
\cite{talviladenjoy}, where the completion is denoted
${\cal A}$.
Thus, if $f\in{\cal A}$ then $f\in{\cal D}'$ (Schwartz
distributions) and
there is a 
function $F\in \tilde{C}$
such that $\langle F',\phi\rangle=-\langle F,\phi'\rangle=
-\int_{-\infty}^\infty F\phi'=\langle f,\phi\rangle$
for all test functions $\phi\in{\cal D}=C^\infty_c(\R)$.
The distributional integral of $f$ is then $\int_a^bf=F(b)-F(a)$
for all $-\infty\leq a\leq b\leq\infty$.
We can
compute the Alexiewicz norm of $f$ via $\|f\|=\sup_{x\in\R}|F(x)|=
\|F\|_\infty$.  If $f\in{\cal D}'$ then $\tau_xf$ is defined by
$\langle\tau_xf,\phi\rangle:=\langle f,\tau_{-x}\phi\rangle
=\langle F',\tau_{-x}\phi\rangle = -\langle F, (\tau_{-x}\phi)'\rangle
=-\langle F, \tau_{-x}\phi'\rangle=-\langle \tau_xF,\phi'\rangle
=\langle (\tau_xF)',\phi\rangle$.  Of course we have $L^1\subset\hk\subset{\cal A}$ and each inclusion is strict.

The theorem only depends on uniform continuity of the primitive
and not on its pointwise differentiability properties so it
also holds in ${\cal A}$.  The same is true for the other theorems
in this paper.
\begin{corollary}
Let $f\in{\cal A}$.  Then $\|\tau_xf-f\|\to 0$ as $x\to 0$.
\end{corollary}

The following theorem gives more precise information on the decay rate 
of $\|\tau_xf-f\|$.

\begin{theorem}\label{theoremdecayrate}
{\rm (a)} Let $\psi\fn(0,1]\to(0,\infty)$ such that $\lim_{x\to 0}\psi(x)=0$.
Then there is $f\in L^1$ such that $\|\tau_xf-f\|\geq\psi(x)$ for all
sufficiently small $x>0$.
{\rm (b)} If $f\in\hk$ and $f\not=0$ a.e. then the most rapid decay is
 $\|\tau_xf-f\|=O(x)$ as $x\to 0$
and this is the best estimate in the sense that if $\tauxf/x\to0$ as
$x\to0$ then $f=0$ a.e.
The implied constant in the order relation is the oscillation of $f$.
\end{theorem}

\bigskip
\noindent
{\bf Proof:} (a) Given $\psi$, define
$\psi_1(x)=\sup_{0<t\leq x}\psi(t)$.  Then $\psi_1\geq\psi$ and
$\psi_1(x)$ decreases to 0 as $x$ decreases to 0.  Define $\psi_2(x)=
\psi_1(1/n)$ when $x\in(1/(n+1),1/n]$ for some $n\in\N$.  Then
$\psi_2\geq \psi$ and $\psi_2$ is a step function that decreases
to 0 as $x$ decreases to 0.  Now let 
$$
\psi_3(x)=
\left[\psi_2(1/(n-1))-\psi_2(1/n)\right]n(n+1)
\left(x-\frac{1}{n+1}\right)+\psi_2(1/n)$$
when $x\in[1/(n+1),1/n]$ for some $n\geq 2$.  Define $\psi_3=\psi_2$ on
$(1/2,1]$.  Then $\psi_3\geq \psi$ and 
$\psi_3$ is a piecewise linear continuous function that decreases
to 0 as $x$ decreases to 0.  Define $f(x)=\psi_3'(x)$ for $x\in(0,1]$
and $f(x)=0$, otherwise.  For $0<x<1$, 
\begin{eqnarray*}
\tauxf & \geq & \left|\int_0^x\left[f(y-x)-f(y)\right]\,dy\right|\\
 & = & \int_0^xf\\
 & = & \psi_3(x)\\
 & \geq & \psi(x).
\end{eqnarray*}
Since $\psi_3$ is absolutely continuous, $f\in L^1$.

(b) Test functions are dense in $\hk$, i.e., for each $f\in\hk$ and 
$\epsilon>0$ there
is $\phi\in\D$ such that $\|f-\phi\|<\epsilon$.  Let $x\in\R$.  Then, since
$\tau_x$ is a linear
isometry,
$\|(\tau_xf-f) -(\tau_x\phi-\phi)\|=\|\tau_x(f-\phi)-(f-\phi)\|<2\epsilon$
and $\|\tau_x\phi-\phi\|-2\epsilon<\tauxf<\|\tau_x\phi-\phi\|+2\epsilon$.  
It therefore suffices to
prove the theorem in $\D$.  Hence, let $f\in\D$ and let $a,b\in\R$.  Write
$F(y)=\int_{-\infty}^yf$.  Then, since $F\in C^2(\R)$, 
\begin{eqnarray*}
\int_a^b(\tau_xf-f) & = & \left[F(b-x)-F(b)\right]-[F(a-x)-F(a)]\\
 & = & -F'(b)\,x+F''(\xi)\,x^2+F'(a)\,x-F''(\eta)\,x^2,
\end{eqnarray*}
for some $\xi,\eta$ in the support of $f$.  Now,
\begin{eqnarray*}
\|\tau_xf-f\| & \geq & \sup_{a,b\in\R}|f(a)-f(b)||x|-2\|f'\|_\infty x^2\\
 & = & {\rm osc}f\, |x|-2\|f'\|_\infty x^2.
\end{eqnarray*}
The oscillation of $f\in\D$ is positive unless $f$ is constant, but there
are no constant functions in $\D$ except 0.  The proof is completed
by noting that $\|\tau_xf-f\|\leq {\rm osc}f\, |x|+2\|f'\|_\infty x^2$
so that $\|\tau_xf-f\|=O(x)$ as $x\to 0$.\Qed

Part (b) is proven in \cite[Proposition~1.2.3]{reiter} for $f\in L^1$.

It is interesting to note that if $f\in \hk$ and $F$ is its primitive then
the function $\tau_xF-F$ is in $\hk$ for each $x\in\R$, 
even though $F$ need not be in $\hk$.

\begin{theorem}
Let $f\in\hk$,  let $F$ be one of its primitives and let $x\in\R$.  
Then the function
$y\mapsto \tau_xF(y)-F(y)$ is in $\hk$ even though none of the primitives of
$f$ need be in $\hk$.  
We have the estimate $\|\tau_xF-F\|\leq \|f\||x|$.  In
general, $\tau_xF-F$ need not be in $L^1$.
However, if $f\in L^1$ then $\tau_xF-F\in L^1$ and 
$\|\tau_xF-F\|_1\leq \|f\|_1|x|$.
\end{theorem}

\bigskip
\noindent
{\bf Proof:}  Let $f\in\hk$ and let $F$ be any primitive.  Since $F$ is
continuous, to prove $\tau_xF-F\in\hk$ we need only show integrability at
infinity.  Let $a,x\in\R$.  Then
\begin{eqnarray*}
\int_0^a(\tau_xF-F) & = & \int_{-x}^{a-x}F-\int_{0}^aF\\
 & = & \int_{-x}^0F-\int_{a-x}^aF\\
 & = &  \int_{-x}^0F -F(\xi)\,x
\end{eqnarray*}
for some $\xi$ between $a-x$ and $a$, due to continuity of $F$.  Hence,
$\lim_{a\to\pm\infty}\int_0^a(\tau_xF-F)= 
\int_{-x}^0F -x\lim_{y\to\pm\infty}F(y)$.
Since $F$ has limits at $\pm\infty$, Hake's theorem shows $\tau_xF-F\in\hk$.
Now let $a,b\in\R$.  Then $\int_a^b(\tau_xF-F)=\int_{a-x}^{a}F-\int_{b-x}^{b}F$.
Since $F$ is continuous, there are $\xi$ between
$a$ and $a-x$ and $\eta$ between $b$ and $b-x$ such that
$\int_a^b(\tau_xF-F)=F(\xi)\,x-F(\eta)\,x=x\int_\eta^\xi f$.  It follows that
$\|\tau_xF-F\|\leq \|f\||x|$.

The example $f=\chi_{[0,1]}$, for which 
$$
F(y)=\int_{-\infty}^y f= \left\{\begin{array}{cl}
0, & y\leq 0\\
y, & 0\leq y\leq 1\\
1, & y\geq 1
\end{array}
\right.
$$
shows that no primitives need not be in $\hk$.
And, if we let $F(y)=\sin(y)/y$, $f=F'$, then for
$x\not=0$,
\begin{eqnarray*}
\tau_xF(y)-F(y) & = & \frac{\sin(y-x)}{y-x}-\frac{\sin(y)}{y}\\
 & \sim & \frac{[\cos(x)-1]\sin(y) -\sin(x)\cos(y)}{y} \text{ as } y\to\infty.
\end{eqnarray*}
Hence, $\tau_xF-F\in\hk\setminus L^1$.

Suppose $f\in L^1$ and $x\geq 0$.  Then,
$|f|\in \hk$ so the theorem gives
$\|\tau_xF-F\|_1\leq \int_{-\infty}^\infty\int_{y-x}^y|f(z)|\,dzdy\leq
\|\, |f|\, \|\,x=\|f\|_1x$.  Similarly, if $x<0$.
\Qed

\begin{example}
{\rm
Let $f$ be $2\pi$-periodic and Henstock--Kurzweil integrable over one
period.  The Poisson
integral of $f$ on the unit circle is 
$$
u(re^{i\theta})=u_r(\theta)=\frac{1-r^2}{2\pi}\int_{-\pi}^\pi\frac{f(\phi)\,d\phi}
{1-2r\cos(\phi-\theta)+r^2}.
$$
Differentiating under the integral sign shows that $u$ is harmonic in
the disc.  And, after interchanging the order of integration, it can be 
seen that
$\|u_r-f\|\to 0$ as $r\to 1^-$.
The Poisson integral defines a harmonic function that takes on the
boundary values $f$ in the Alexiewicz norm.  For details see 
\cite{talvilapoisson}.\Qed
}
\end{example}

Now we consider continuity in weighted Alexiewicz norms.  First we need
the following lemma.  Lebesgue measure is denoted $\lambda$.

\begin{lemma}\label{lemma1}
For each $n\in\N$, suppose $g_n\fn\R\to\R$ and $g_n\chi_E\to g\chi_E$
in measure for some set $E\subset \R$ of positive measure and function $g$
of bounded variation.  If
$Vg_n\leq M$ for all $n$ then $g_n$ is uniformly bounded on $\R$.
\end{lemma}

\bigskip
\noindent
{\bf Proof:} Define $S_n=\{x\in E\mid |g_n(x)-g(x)|>1\}$.
Then $\lambda(S_n)\to 0$ as $n\to \infty$.  There is $N\in\N$ such
that whenever $n\geq N$ we have $\lambda(E\setminus S_n)>0$.
Since $g\in\bv$,  $g$ is bounded. 
Let $n\geq N$.  There is $x_n\in E\setminus S_n$ such that
$|g(x_n)|\leq \|g\|_\infty$.  Therefore, $|g_n(x_n)|\leq 1+\|g\|_\infty$.
Let $x\in\R$.  Then 
$|g_n(x)-g_n(x_n)|\leq Vg_n\leq M$.
So, $|g_n(x)|\leq M+1+\|g\|_\infty$.  Hence, $\{g_n\}$ is uniformly
bounded.\Qed

\begin{theorem}\label{theoremweighted}
Let $w\fn\R\to(0,\infty)$.  Define $g_x\fn\R\to(0,\infty)$
by $g_x(y)=w(y+x)/w(y)$ for each $x\in\R$.
Then $\|(\tau_xf-f)w\|\to 0$ as $x\to 0$ 
for all $f\fn\R\to\R$ such that $fw\in\hk$ if and only if
$g_x$ is essentially bounded and of essential bounded
variation, uniformly as $x\to 0$, and $g_x\to 1$ in measure
on compact intervals as $x\to0$.
\end{theorem}

\noindent
{\bf Proof:} Let $G(x)=\int_{-\infty}^xfw$. Let $x,\alpha,\beta\in\R$.  Then
\begin{eqnarray*}
\lefteqn{
\int_\alpha^\beta\left[f(y-x)-f(y)\right]w(y)\,dy    }\\
& = & \int_{\alpha-x}^{\beta-x}f(y)w(y)\,dy - \int_{\alpha}^{\beta}f(y)w(y)\,dy
+ \int_{\alpha-x}^{\beta-x}f(y)\left[w(y+x)-w(y)\right]\,dy\\
 & = & \left[G(\beta-x)-G(\beta)\right] -
\left[G(\alpha-x)-G(\alpha)\right]+ \int_{\alpha-x}^{\beta-x}
f(y)w(y)\left[g_x(y)-1\right]\,dy.
\end{eqnarray*}
Since $G$ is uniformly continuous on $\R$, we have
$\|(\tau_xf-f)w\|\to 0$ if and only if the supremum of $|\int_{a}^{b}
f(y)w(y)[g_x(y)-1]\,dy|$ over $a,b\in\R$ has limit $0$ as $x\to 0$,
i.e., $\|fw(g_x-1)\|\to 0$.  Given $h\in\hk$ we can always take
$f=h/w$.  Hence,  the theorem now follows from Lemma~\ref{lemma1} (easily
modified for the case of essential boundedness and essential variation) and
the necessary
and sufficient condition for convergence in norm given in 
\cite[Theorem~6]{mohantytalvila}.
\Qed

\begin{corollary}
Suppose that for each compact interval $I$ there are real numbers
$0<m_I<M_I$ such that $m_I<\|w\|_\infty<M_I$; $w$ is continuous in measure
on $I$; $w\in\bv_{loc}$.  Then for all $f\fn\R\to\R$ such that $fw\in\hk$
we have $\|(\tau_xf-f)w\|\to 0$ as $x\to 0$.
\end{corollary}

\bigskip
\noindent
{\bf Proof:} Fix $\epsilon>0$.  Let $I$ be a compact interval for which
$0<m_I<\|w\|_\infty<M_I$.  Define 
\begin{eqnarray*}
S_x & :=  & \{y\in I\mid |g_x(y)-1|>\epsilon\}\\
 & = &  \{y\in I\mid |w(y+x)-w(y)|>\epsilon\, w(y)\}\\
 & \subset &  \{y\in I\mid |w(y+x)-w(y)|>\epsilon\,m_I\}
\quad\text{except for a null set}.
\end{eqnarray*}
Since $w$ is continuous in measure on $I$ we have $\lambda(S_x)\to 0$
as $x\to 0$ and $g_x\to 1$ in measure on $I$.

Using 
$$
g_x(s_n)-g_x(t_n)=\frac{w(s_n+x)-w(t_n+x)}{w(t_n)}
-\frac{w(s_n+x)[w(s_n)-w(t_n)]}{w(s_n)w(t_n)}
$$
we see that $V_Ig_x\leq V_{I+x}w/m_I+M_IV_Iw/m^2_I$ where
$I+x=\{y+x\mid y\in I\}$ and $V_{I}w$ is the variation of $w$
over interval $I$.  Hence, $g_x$ is of uniform bounded variation
on $I$.  With Lemma~\ref{lemma1} this then gives the hypotheses of the theorem. \Qed

Of course, we are allowing $w$ to be changed on a set of measure $0$ so
that $w$ is of bounded variation rather than just equivalent to a function
of bounded variation.  This redundancy can be removed by replacing $w$
with its limit from the right at each point so that $w$ is right continuous.

As pointed out in \cite{mohantytalvila}, convergence of $g_x$ to $1$ in measure
on compact intervals in the theorem 
can be replaced by convergence in $L^1$ norm:
For each compact interval $I$, $\|(g_x-1)\chi_I\|_1\to 0$ as $x\to 0$. 
In the  corollary, we can replace continuity in measure with the condition:
As $x\to 0$, $\int_I|\tau_xw-w|\to 0$ for each compact interval $I$.

The first two conditions in the corollary are necessary.  Suppose
$\|(\tau_xf-f)w\|\to 0$ whenever $fw\in\hk$.  Then the essential
infimum of $w$ must be positive on each compact interval.  If there
is a sequence $a_n\to a\in\R$ for which $w(a_n)\to 0$ then
${\rm esssup}_{|x|<\delta} \sup_{n}g_x(a_n)=\infty$ for each $\delta>0$
unless $w=0$ a.e. in a neighbourhood of $a$.  Similarly, the essential
supremum of $w$ must be finite on compact intervals.  This asserts the
existence of $m_I$ and $M_I$ in the corollary.
Also, let $I$ be a compact interval on which $0<m<\|w\|_\infty<M$.  Let
$\epsilon>0$ and define
\begin{eqnarray*}
T_x & :=  & \{y\in I\mid |w(y+x)-w(y)|>\epsilon\}\\
 & \subset &  \{y\in I\mid |g_x(y)-1|>\epsilon/M\}
\quad\text{except for a null set}.
\end{eqnarray*}
Hence, $w$ is continuous in measure.

It is not known if $\|(\tau_xf-f)w\|\to 0$ for all $f$ such
that $fw\in\hk$ implies $w\in\bv_{loc}$.  The example
$w(y)=e^y$ shows $W$ need not be of bounded  variation and can
have its infimum zero and its supremum infinity.  For, note
that $g_x(y)=\exp(y+x)/\exp(y)=e^{x}$ and so satisfies the
conditions of the theorem.  And, by the corollary,
$w(y)=1$ for $y<0$ and $w(y)=2$ for $y\geq 0$ is a valid weight 
function.  Hence, $w$ need not be continuous.

\begin{example}
{\rm 
Let $w(y)=1/(y^2+1)$.  A calculation shows that the variation of
$y\mapsto w(y+x)/w(y)$ is $2|x|\sqrt{x^2+1}$ so $w$ is a valid weight for 
Theorem~\ref{theoremweighted}. The
half-plane Poisson kernel is $\Phi_y(x)=w(x/y)/(\pi y)$.  For $f\fn\R\to\R$ the
Poisson integral of $f$ is $u_y(x)= (\Phi_y\ast f)(x)=
\frac{y}{\pi}\int_{-\infty}^\infty\frac{f(t)\,dt}{(x-t)^2+y^2}$.  Define
$\Psi_z(t)=\Phi_y(x-t)/w(t)$ for $z=x+iy$ in the upper half-plane, i.e.,
$x\in\R$ and $y>0$.  For fixed $z$ both $\Psi_z$ and $1/\Psi_z$ are of
bounded variation on $\R$.  Hence, necessary and sufficient for the existence
of the Poisson integral on the upper half-plane is $fw\in\hk$.

Define $G(t)=\int_{-\infty}^tfw$.  Integrate by parts to get 
$u_y(x)=y\,G(\infty)/\pi -\int_{-\infty}^\infty G(t)\Psi_z'(t)\,dt$.  Since $G$ is
continuous on the extended real line (with $G(\infty):=\lim_{t\to\infty}G(t)$),
dominated convergence now allows differentiation under the integral.  This
shows $u_y(x)$ is harmonic in the upper half-plane.

Neither $f$ nor $u_y$ need be in $\hk$. For example, the Poisson integral
of $1$ is $1$.
But, we have the boundary values taken on in the weighted norm:
$\|(u_y-f)w\|\to 0$ as $y\to 0^+$.  We sketch out the proof, leaving the
technical detail of interchanging repeated integrals for publication elsewhere.
For $a,b\in\R$ we then
have
\begin{eqnarray*}
\int_a^b\left[u_y(t)-f(t)\right]w(t)\,dt & = &
\int_a^b\left\{(f\ast\Phi_y)(t) - f(t)\int_{-\infty}^\infty 
\Phi_y(s)\,ds\right\}w(t)\,dt\\
 & = & \int_{-\infty}^\infty \Phi_y(s)\int_a^b\left[f(t-s)-f(t)\right]w(t)\,dt\,ds.
\end{eqnarray*}
Therefore, 
$\|(u_y-f)w\|\leq \int_{-\infty}^\infty\Phi_y(s)\|(\tau_sf-f)w\|\,ds$.
But, $s\mapsto \|(\tau_sf-f)w\|$ is continuous at $s=0$.  By the usual
properties of the Poisson kernel (an approximate identity), we have
$\|(\tau_sf-f)w\|\to 0$ as $y\to 0^+$.
\Qed
}
\end{example}

\end{document}